%% file: main.tex
\documentclass{article}

\usepackage{lmodern}

\usepackage{graphicx} 
\graphicspath{{images/}{}}
\usepackage{float}
\usepackage{longtable}

\usepackage{caption}
\usepackage{subcaption}

\usepackage{indentfirst}

\usepackage{amsmath}
\usepackage{amssymb}
\usepackage{mathtools}
\usepackage{amsthm}

\usepackage{verbatim}

\usepackage[colorlinks=true, allcolors=blue]{hyperref}

\usepackage[style=numeric-comp]{biblatex}
\usepackage{soul} 

\theoremstyle{conjecture}
\newtheorem*{conjecture}{Conjecture}

\title{Porosity and topological properties of triply periodic minimal surfaces}

\makeatletter
\def\@maketitle{%
  \newpage
  \null
  \vskip 2em%
  \begin{center}%
  \let \footnote \thanks
    {\Large\bfseries \@title \par}%
    \vskip 1.5em%
    {\normalsize
      \lineskip .5em%
      \begin{tabular}[t]{c}%
        \@author
      \end{tabular}\par}%
    \vskip 1em%
    {\normalsize \@date}%
  \end{center}%
  \par
  \vskip 1.5em}
\makeatother

\usepackage{authblk}
\author[1]{Pavel Snopov}
\author[2]{Sergei Ermolenko}
\affil[1]{Institute for Information Transmission Problems, RAS}
\affil[2]{Voronezh State University} 
\date{}

\begin{document}

\maketitle

\begin{abstract}
\noindent Triple periodic minimal surfaces (TPMS) have garnered significant interest due to their structural efficiency and controllable geometry, making them suitable for a wide range of applications. This paper investigates the relationships between porosity and persistence entropy with the shape factor of TPMS. We propose conjectures suggesting that these relationships are polynomial in nature, derived through the application of machine learning techniques. This study exemplifies the integration of machine learning methodologies in pure mathematical research. Besides the conjectures, we provide the mathematical models that might have the potential implications for the design and modeling of TPMS structures in various practical applications.

\end{abstract}

\input{introduction}

\input{math_models}

\input{persistence}

\input{calc_porosity}

\input{ph_of_tpms}

\input{regression_part}

\input{conclusion}

\printbibliography[title={Bibliography}]

\end{document}

%% file: introduction.tex
\section{Introduction}

Triple periodic minimal surfaces (TPMS) are cellular structures with geometry that can be tailored for desired mechanical responses, offering structural efficiency by minimizing surface area and optimizing mass distribution \cite{Asb21, Asb22}. Despite traditionally viewing pores as defects, natural porous materials exhibit strong mechanical and mass transfer properties, making TPMS particularly valuable for modeling and manufacturing.

In medicine, TPMS-designed implants, such as porous meniscus implants, can reduce compressive and shear stresses compared to solid implants \cite{Asb23}. The physical properties of these models are influenced by their porosity, determined during mathematical and computer modeling. However, manufacturing such complex topologies remains challenging \cite{Asb220}. Traditional fabrication methods struggle to replicate natural porous structures accurately \cite{Asb213}, and ensuring reliable manufacturing quality and broad applicability of these materials is an ongoing challenge \cite{Asb214}.

Persistent homology \cite{comptop, persistence} is a novel tool in applied topology, integral to topological data analysis, utilizing algebraic topology to analyze data. It provides a stable method to capture changes in algebraic invariants from topological objects derived from data.

In this paper, we analyze geometric and topological properties of two specific TPMS, Schwarz Primitive and Gyroid. We calculate their porosity, compute persistence diagrams, and examine the relationship between porosity and 1-persistence entropy with the shape factor of TPMS Our findings suggest a potential polynomial dependence between these properties and the shape factor, leading us to propose related conjectures. Finally, we are publicly releasing the code\footnote{https://github.com/Snopoff/Persistent-Homology-of-Triple-Periodic-Minimal-Surfaces}, hoping to faciliate the applicability of the TDA methods to the study of TPMS and other porous structures.

%% file: math_models.tex
\section{Preliminaries} 
\subsection{Mathematical model of minimal surfaces}
\label{100}

Triple periodic minimal surfaces (TPMS) are a prominent part of the well-established minimal surface problem, the Plateau problem in three-dimensional Euclidean space ($\mathbb{R}^3$). The core of this problem lies in investigating the shape and configuration of a soap film when constrained by a wire frame of a given geometry. Consequently, the described minimal surface represents a surface whose mean curvature is identically zero.

To model TPMS, a mathematical framework is employed that allows for the representation of a surface with a specified triple periodic symmetry and minimal area. The modeling of triple periodic minimal surfaces frequently relies on the minimal surface principle, which stipulates that the surface with the smallest area is in equilibrium and possesses the least energy. Various mathematical techniques, such as variational methods \cite{Asb222} and differential equation approaches \cite{Asb26}, are utilized to determine minimal surfaces.

Poisson established that if a two-dimensional smooth surface \(M^2\) is embedded in \(R^3\) and represents the interface between two equilibrium media, then the mean curvature of the surface is constant and equal to:

$$H = \frac{1}{2}(k_1 + k_2)$$

where \(H\) denotes the mean surface curvature, and $k_1$ and $k_2$ are the principal curvatures of the surface at a given point. The condition of minimal area leads to the requirement that \(H\)=0. Consequently, a surface with a mean curvature of zero is considered minimal \cite{Asb24}.

In the case of a soap film:

$$H = \frac{p}{2\sigma} = 0$$

where $p$ represents the pressure difference across the surface, $\sigma$ - surface tension. When the pressure on both sides of the film is equal, i.e., $P_1 = P_2$, then $p$ and \(H\) is equal zero. This is because if there is no pressure differential on either side of the surface, there is no driving force for the film to bend or deform \cite{Asb223}.

The mean curvature formula allows for the calculation of the curvature of the surface at each point. This property of mean curvature makes the above formula a defining characteristic of minimal surfaces \cite{Asb25}.

Triple periodic symmetry implies that the surface must be symmetrical about three axes (such as the \(x\), \(y\), and \(z\) directions) and repeat three times in each of those directions.

Assuming the surface is repeated three times in each direction, the boundary conditions must account for the fact that the values of the surface parameters must be the same when shifted by an integer number of periods. For the specific case of triple periodic symmetry, the equations of the boundary conditions can be written.

Let \(f(x, y, z)\) be a function describing the surface. Then, if the surface is periodic in the x direction with period a, the boundary condition can be written as:
\[f(x, y, z) = f(x+a, y, z)\]

Similar boundary conditions can be expressed for the periodicity in the \(y\) and \(z\) directions.

The minimal surface area is typically calculated using an integral formula that accounts for the elementary surface energy and geometric derivatives of the surface parameters.

Since the minimal surface equations are often intricate and lack analytical solutions, the utilization of numerical methods for approximating solutions, such as the finite element method, is required. This approach allows for the division of a surface into small elementary regions called finite elements and the approximation of the surface of each element using local functions \cite{Asb218, Asb216}.

In this study, the following triply periodic minimal surfaces (TPMS) were investigated: Gyroid and Schwarz Primitive. As mentioned previously, TPMS are expressed by a mathematical formula that describes their three-dimensional geometry. This fact is crucial when modeling objects using TPMS, as by altering the parameters of a particular mathematical model, one can obtain the form necessary for a specific application. The first examples of TPMS are surfaces described by Schwartz in 1865 with the following equation:
\[f(x, y, z) = \cos{ax} + \cos{by} + \cos{cz} + d\]
As for Gyroid, its equation has a more complex form:
\[f(x, y, z) = \sin{ax}\cos{by} + \sin{by}\cos{cz} + \sin{cz}\cos{ax} + d\]
where a, b, c are parameters affecting the characteristic size of the unit cell of a periodic structure; d - parameter influencing the volume fraction of the space cut off by TPMS. In TPMS models, the dimensional parameters of the unit cell of a periodic structure influence the axial periodicity of the surface. This is an essential feature when developing physical models for applied problems. TPMS possess several advantages due to the previously mentioned properties and characteristics.

%% file: persistence.tex
\subsection{Persistence}
\label{subsec:persistence} \index{Persistent Homology}
Recall that a simplicial complex $X$ is a collection of $p$-dimensional simplices, i.e., vertices, edges, triangles, tetrahedrons, and so on. Simplicial complexes can be seen as a generalization of graphs, which are essentially pairs of vertices ($0$-simplices) and edges ($1$-simplices) that can represent higher-order interactions. A family of increasing simplicial complexes

\[
\emptyset \subseteq X_0 \subseteq X_1 \subseteq \cdots \subseteq X_{n-1} \subseteq X_n
\]

is called a \textit{filtration}.

The idea behind \textit{persistence} \cite{persistence} is to track the evolution of simplicial complexes over the filtration. {\it Persistent homology} allows one to trace the changes in homology vector spaces\footnote{Usually, homology groups are considered with integral coefficients, but in the realm of persistence, homology groups are taken with coefficients in some field, for example, $\mathbb{Z}_p$ for some large prime number $p$. Hence, the homology groups become homology vector spaces.} of simplicial complexes that are present in the filtration. Given a filtration $\{X_t\}_{t=0}^n$, the homology functor $H_p$ applied to the filtration produces a sequence of vector spaces $H_p(X_t)$ and maps $i_*$ between them

\[
H_p(X_*): H_p(X_0) \xrightarrow{i_0} H_p(X_1) \xrightarrow{i_1} \cdots \xrightarrow{i_{n-1}} H_p(X_n).
\]

Each of these vector spaces encodes some information about the simplicial complex $X$ and its subcomplexes $X_i$. For example, $H_0$ generally encodes the connectivity of the space (or, in data terms, $H_0$ encodes the clusters of data), $H_1$ encodes the presence of 1-cycles, i.e., loops, $H_2$ encodes the presence of 2-cycles, and so on. Persistence allows one to track the generators of each vector space via the induced maps. Some generators will vanish, while others will persist. Those that persist are likely the most important, as they represent features that truly exist in $X$. Therefore, persistent homology allows one to gain information about the underlying topological space via the sequence of its subspaces, the filtration.

In algebraic terms, the sequence of vector spaces $H_p(X_t)$ and maps $i_*$ between them can be seen as a representation of a quiver $I_n$. From the representation theory of quivers \cite{representation_theory}, it is known, due to Gabriel, that any such representation is isomorphic to a direct sum of indecomposable interval representations $I[b_i, d_i]$, that is,

\[
H_p(X_*) \simeq \bigoplus_i I[b_i, d_i].
\]

\begin{figure}
\centering
\begin{subfigure}{0.5\textwidth}
  \centering
  \includegraphics[width=0.9\linewidth]{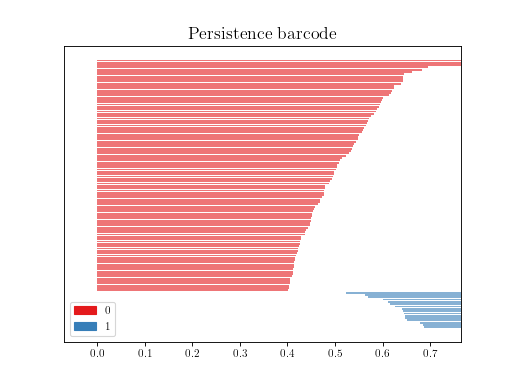}
  \caption{An example of a barcode}
  \label{fig:barcode}
\end{subfigure}%
\begin{subfigure}{0.5\textwidth}
  \centering
  \includegraphics[width=0.9\linewidth]{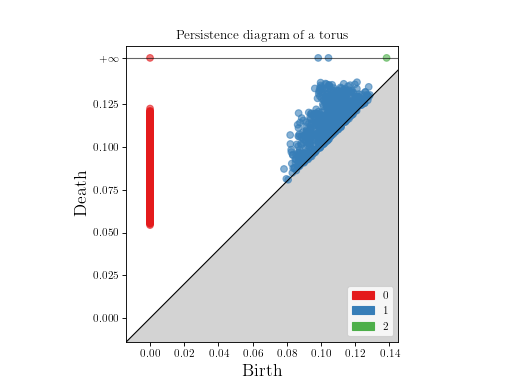}
  \caption{An example of a persistence diagram}
  \label{fig:diagram}
\end{subfigure}
\caption{The examples of the persistent homology visualizations \cite{gudhi}}
\label{fig:persistence}
\end{figure}

The pairs $(b_i, d_i)$ represent the persistence of topological features, where $b_i$ denotes the time of birth and $d_i$ denotes the time of death of the feature. These pairs can be visualized via {\it barcodes} (see Fig. \ref{fig:barcode}), where each bar starts at $b_i$ and ends at $d_i$.

This information is also commonly represented using {\it persistence diagrams} \cite{comptop} (see Fig. \ref{fig:diagram}). A persistence diagram is a (multi)set of points in the extended plane $\overline{\mathbb{R}^2}$, which reflects the structure of persistent homology. Given a set of pairs $(b_i, d_i)$, each pair can be considered as a point in the diagram with coordinates $(b_i, d_i)$. Thus, a persistence diagram is defined as

\[
\mathrm{dgm}(H_p(X_*)) = \{ (b_i, d_i) : I[b_i, d_i] \text{ is a direct summand in } H_p(X_*) \}.
\]

\subsubsection{Filtration Methods}
\label{subsec:filtration} \index{Filtration Methods}
There are various point cloud filtration methods in topological data analysis. A commonly used method is the {\it Vietoris-Rips} filtration $VR_r(X)$ \cite{tda4ds}. This filtration assigns a simplicial complex for each non-negative $r$, where $k+1$ points form a $k$-simplex if the pairwise distances between these points are less than $r$.
Alternatively, one can visualize this process by considering closed balls of radius $r$ centered at each point in the point cloud. When $k+1$ balls intersect, a $k$-simplex is formed in the complex.

Another approach involves constructing a filtration using {\it Alpha complexes} \cite{comptop}. An Alpha complex is a subcomplex of the Delaunay triangulation, obtained by intersecting the closed balls, as in the Vietoris-Rips complex, with the Voronoi cells. This type of filtration tends to produce complexes with significantly fewer simplices, thereby accelerating the persistent homology calculations.

\subsubsection{Persistence Entropy}
\label{subsec:entropy} \index{Persistence Entropy}
Persistence diagrams offer a detailed visual representation of the topology of point clouds. However, incorporating them into machine learning models poses significant challenges due to their structure. To utilize the information contained in persistence diagrams within predictive models, it is essential to transform the data into a suitable representation, typically a vector within a vector space. This transformation process is referred to as {\it vectorization}.

One approach to vectorize a persistence diagram involves mapping it to a real number. {\it Persistence entropy} \cite{entropy} is a method for achieving this. Persistence entropy is a specialized form of Shannon entropy and is calculated as follows:

\[
PE_k(X) \coloneqq -\sum_{(b_i, d_i) \in D_k} p_i \log(p_i),
\]

where

\[p_i = \frac{d_i - b_i}{\sum_{(b_i, d_i) \in D_k} (d_i - b_i)} \quad \text{and} \quad D_k \coloneqq \mathrm{dgm}(H_k(X_\bullet)).\]

This numerical characteristic of a persistence diagram has several advantageous properties. Notably, under certain mild assumptions, it is stable. Stability implies the existence of a bound that <<regulates>> the perturbations caused by noise in the input data \cite{entropy_stability}.


%% file: calc_porosity.tex
\section{Porosity calculations of TPMS models with different shape factor}\label{103}

Let us consider the problem of estimating the dependence of the porosity of the TPMS model on changes in the free term, which will be further denoted by us as parameter d and referred to as the shape factor. It is assumed that this parameter influences the volume fraction of the space cut off by the surface. We will use the models from section \ref{100} as the object of investigation.

At the initial stage, we will demonstrate whether there exists a linear or other mathematically interpretable dependence of the value of the parameter d on the porosity of the Schwarz and Gyroid surfaces. First, using the Python programming language, we will initialize the random values of this parameter utilizing the built-in random number module of the Numpy library.

To obtain a point cloud of the studied TPMS surfaces and calculate their porosity, we will employ the PyVista module \cite{pyvista}. This module can be used to create scientific visualizations for presentations and research publications, and as a helper module for other Python modules that depend on 3D mesh rendering. Let's generate points along the three axes x, y, z and obtain a cloud of points. 1000 points along each axis were used as input.

Next, we will substitute the obtained values into the equations of the TPMS under investigation and acquire a three-dimensional cloud of points in the Numpy array. Using the PyVista module's volume function, we will calculate the volume of the surface. Then, we will find the volumetric porosity of the surface as follows:

\[Vp = \frac{Vc - Vs}{Vc}\]

where \(Vc\) represents the volume of the entire space (a cube) occupied by the surface,
\(Vs\) is the calculated model volume of the TPMS

We will store the found porosity value and the arrays containing the point clouds in the corresponding Python lists. This way, we will calculate the volumetric porosity for 20 randomly selected free terms (the parameter d) using two models: the Schwarz and Gyroid surfaces.

The results of the calculations are presented in Table 1

\begin{longtable}{|c|c|c|c|}
    \caption{Porosity calculations results}
    \label{tab:porosity} \\
    \hline
        № & Shape factor & Porosity of Schwarz & Porosity of Gyroid\\ \hline 
        \endfirsthead
        
        \multicolumn{4}{c}%
        {\tablename\ \thetable\ -- \textit{Continued from previous page}} \\
        \hline
        1 & 2 & 3 & 4\\ \hline
        \endhead
        
        \hline
        \multicolumn{4}{r}{\textit{Continued on next page}} \\
        \endfoot
        
        \hline
        \endlastfoot
        
        1 & -0.985086 & 0.796602 & 0.806768 \\ \hline
2 & -0.834337 & 0.777147 & 0.831175 \\ \hline
3 & -0.602925 & 0.759682 & 0.859706 \\ \hline
4 & -0.357099 & 0.730533 & 0.887477 \\ \hline
5 & -0.350885 & 0.730399 & 0.887931 \\ \hline
6 & -0.262321 & 0.707480 & 0.896708 \\ \hline
7 & -0.086151 & 0.681178 & 0.995867 \\ \hline
8 & 0.184648 & 0.667845 & 0.989362 \\ \hline
9 & 0.429486 & 0.689464 & 0.989704 \\ \hline
10 & 0.446178 & 0.685388 & 0.992671 \\ \hline
11 & 0.465395 & 0.690117 & 0.991798 \\ \hline
12 & 0.583936 & 0.671619 & 0.984828 \\ \hline
13 & 0.642550 & 0.664104 & 0.975655 \\ \hline
14 & 0.713229 & 0.670550 & 0.974905 \\ \hline
15 & 0.721315 & 0.666064 & 0.980358 \\ \hline
16 & 0.778426 & 0.673395 & 0.975438 \\ \hline
17 & 0.867778 & 0.680311 & 0.979937 \\ \hline
18 & 0.873779 & 0.671111 & 0.980480 \\ \hline
19 & 0.953246 & 0.686131 & 0.995606 \\ \hline
20 & 0.959561 & 0.686817 & 0.996571 \\ \hline
\end{longtable}

Graphs of the obtained dependence of porosity on the parameter d, that is, the shape factor for the Schwarz and Gyroid models are presented below.

\begin{figure}[H]
    \centering
    \includegraphics[width=1\linewidth]{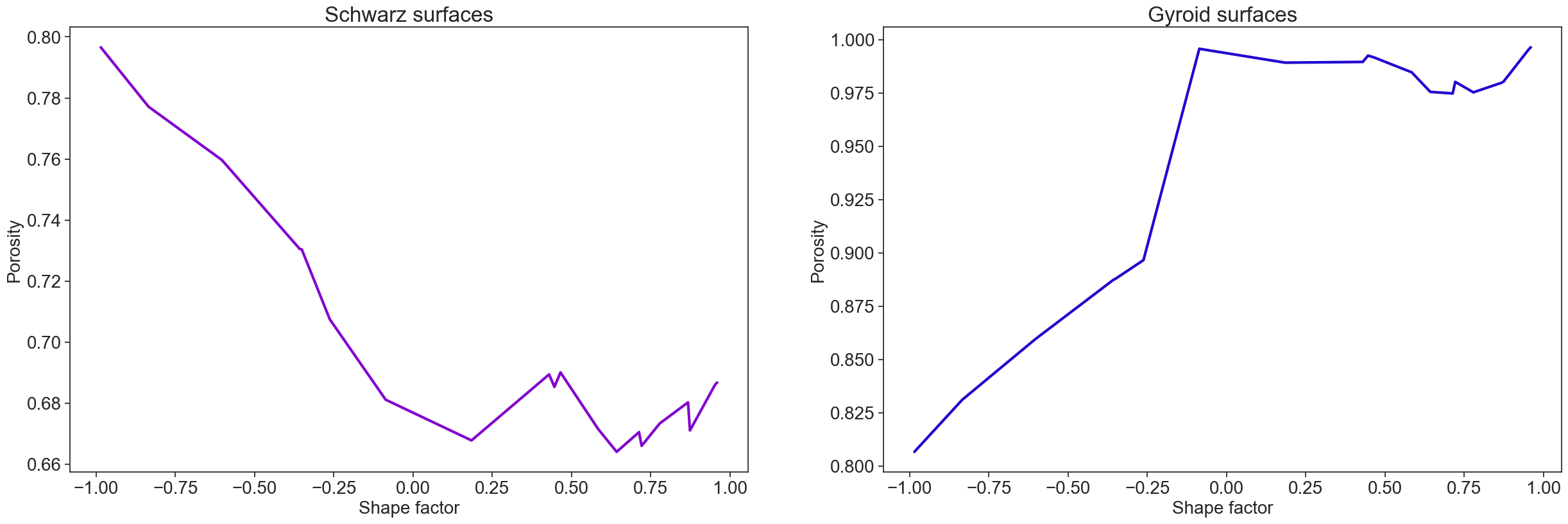}
    \caption{Dependence of TPMS porosity on shape factor}
    \label{fig:enter-label}
\end{figure}

As can be observed from Figure \ref{fig:enter-label}, the porosity of the Schwarz surface decreases with increasing values of parameter d. Moreover, this trend is more pronounced for negative values of the parameter, while for positive values, the porosity indicator tends to stabilize. Regarding the Gyroid surface, an increase in porosity is observed with increasing parameter values, reaching near-maximum levels. Further, a slight decrease is observed for positive parameter values, but the porosity generally remains at approximately the same level. In the investigated models, an obvious linear relationship between porosity and parameter d was not discernible.

\begin{figure}[H]
\centering
\begin{subfigure}{.5\textwidth}
  \centering
  \includegraphics[width=.9\linewidth]{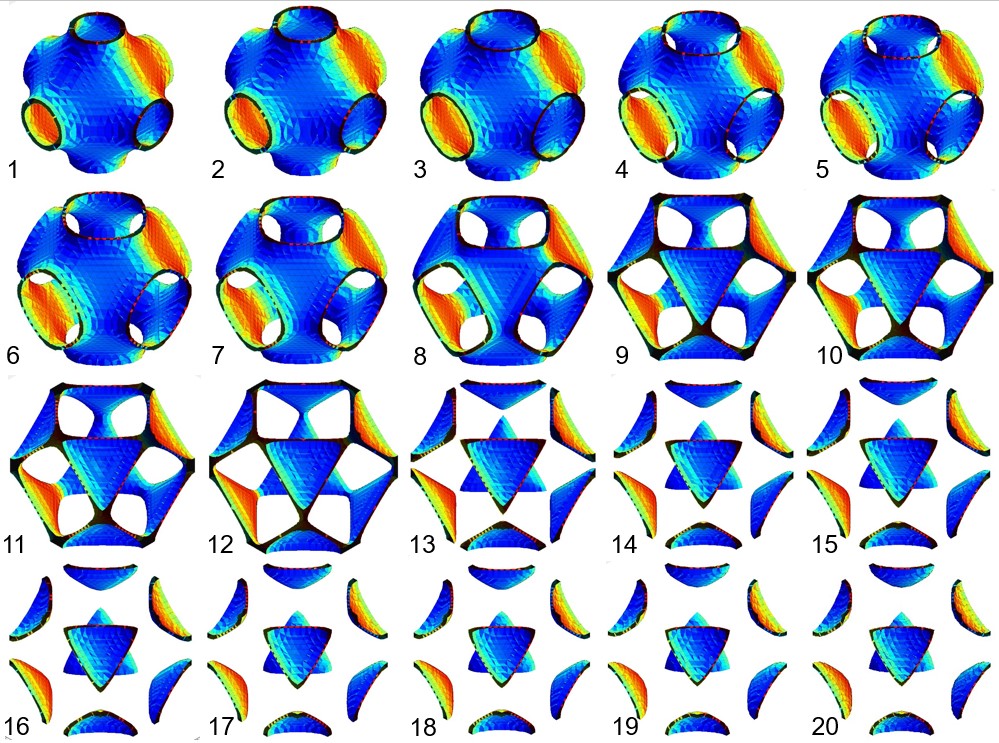}
  \caption{Schwarz surfaces}
  \label{fig:schwarz}
\end{subfigure}%
\begin{subfigure}{.5\textwidth}
  \centering
  \includegraphics[width=.95\linewidth]{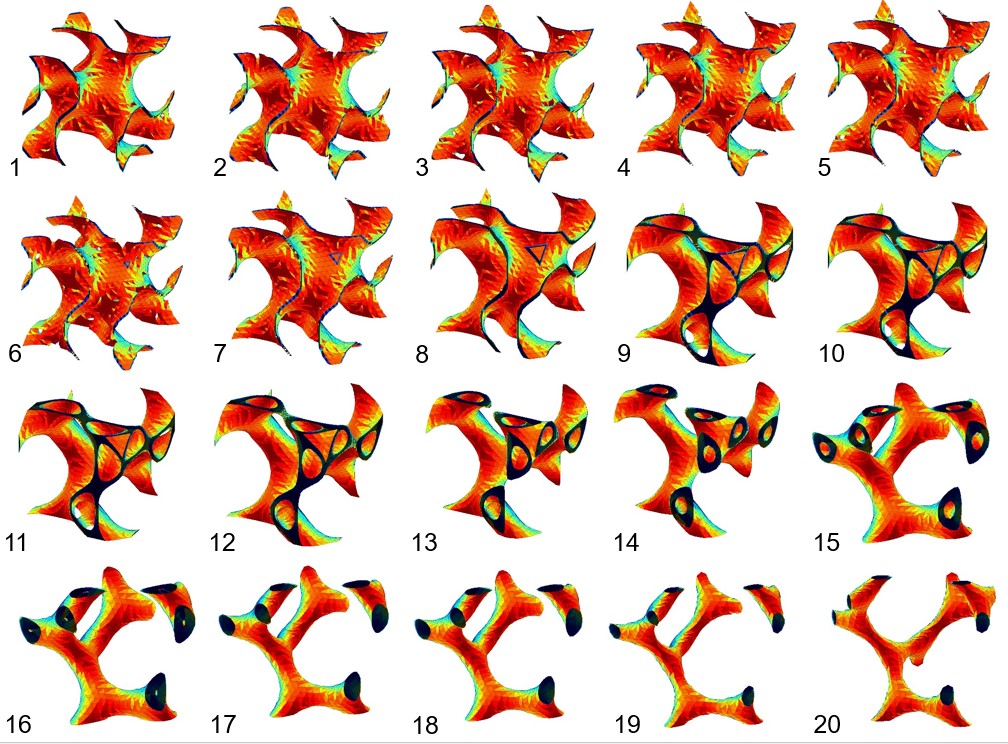}
  \caption{Gyroid surfaces}
  \label{fig:gyroid}
\end{subfigure}
\caption{TPMS surfaces: the number at the bottom-left corner corresponds to the number in the Table \ref{tab:porosity}.}
\label{fig:test}
\end{figure}

We also generated a mesh of the resulting surfaces for each d value using the TPMS-Designer software, an open-source MATLAB toolkit for creating, analyzing, and visualizing TPMS-like structures and other 3D objects \cite{Asb219}. The resulting surface images are presented in Figures \ref{fig:schwarz} and \ref{fig:gyroid}.

%% file: ph_of_tpms.tex
\subsection{Topology estimation of TPMS}
\label{subsec:ph_of_tpms} \index{Persistent Homology of TPMS}
Let's attack the problem of estimating the changes in topology varying the shape factor. In order to do this, we compute the persistent homology of different triple periodic minimal surfaces. The alpha complexes were used in order to compute persistent homology. We used {\it Gudhi} package \cite{gudhi} for persistent homology calculations of point cloud approximations of the surfaces.

\subsubsection{Schwarz}
\begin{figure}[H]
\centering
\begin{subfigure}{\textwidth}
  \centering
  \includegraphics[width=\textwidth]{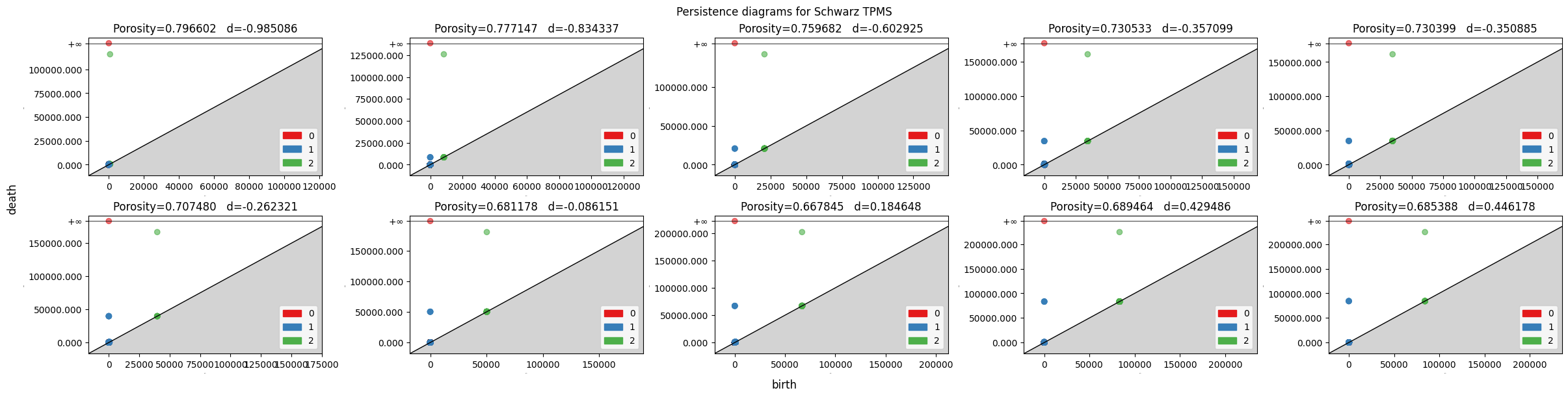}
\end{subfigure}

\vspace{3ex}

\begin{subfigure}{\textwidth}
  \centering
  \includegraphics[width=\textwidth]{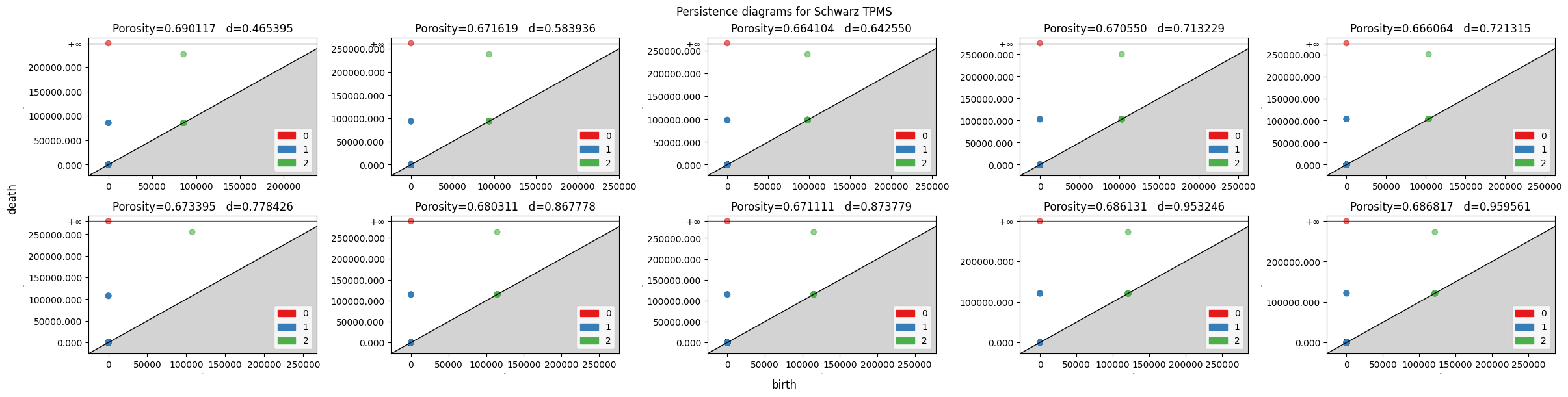}
\end{subfigure}
\caption{Persistence diagrams for Schwarz surfaces}
\label{fig:pd_schwarz}
\end{figure}

The results of the persistent homology calculations are shown in \ref{fig:pd_schwarz}. Notice how topology changes when shape factor changes. In the case of the Schwarz surface, the birgth time of persistent homology generators in 2-nd grading shifts together with the death time. Also notice that besides the 2-persistent homology, the 1-persistent homology generator death time increases as the shape factor increases.

\subsubsection{Gyroid}
\begin{figure}[H]
\centering
\begin{subfigure}{\textwidth}
  \centering
  \includegraphics[width=\textwidth]{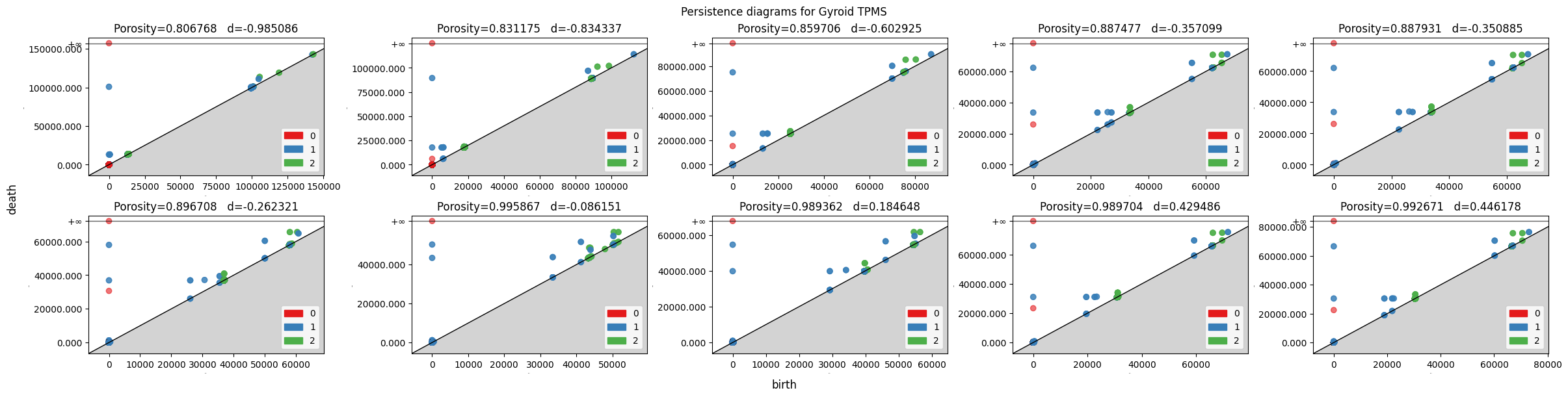}
\end{subfigure}

\vspace{3ex}

\begin{subfigure}{\textwidth}
  \centering
  \includegraphics[width=\textwidth]{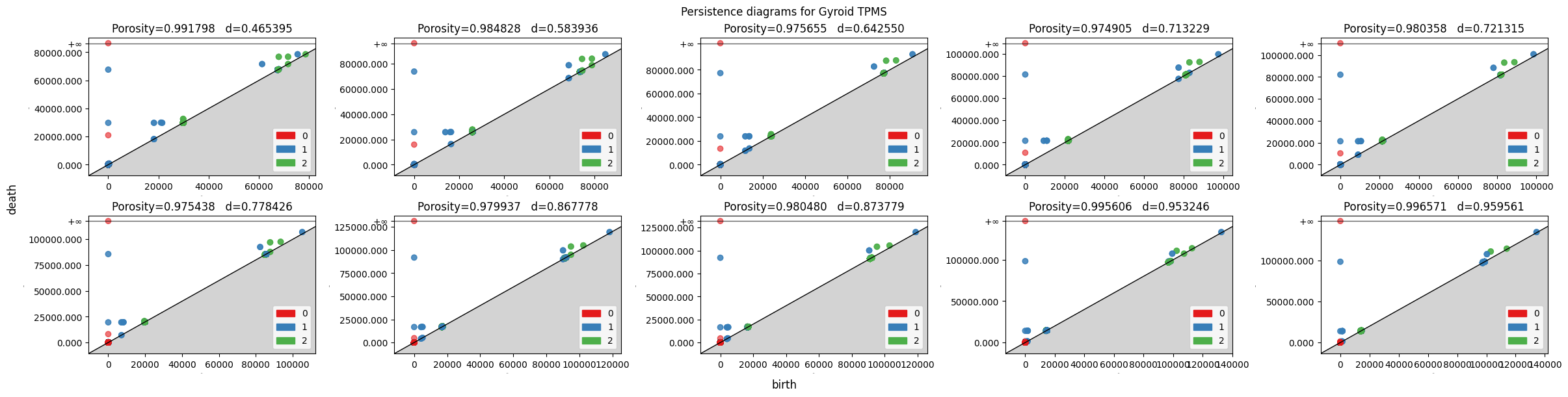}
\end{subfigure}
\caption{Persistence diagrams for Gyroid surfaces}
\label{fig:pd_gyroid}
\end{figure}

The results of the persistent homology calculations are shown in \ref{fig:pd_gyroid}. The picture is different from the Schwarz. When shape factor increases from $-1$ to $0$, the number of generators in all persistent homology grading increases. This means that when shape factor increases, the <<complexity>> of the surface increases as well. But when the shape value increases too much, from $0$ to $1$, the surface' complexity decreases, and the persistence diagram when shape factor equals $1$ is similar to the persistence diagram when shape factor is $0$. This leads to the hypothesis of some quadratic dependency between the shape factor and persistence, which is uncovered in the \ref{subsec:regression_part}.

%% file: regression_part.tex
\section{Exploring Relationships in TPMS} \label{subsec:regression_part} \index{Exploring Relationships in TPMS}

Now we study the relationship between the shape factor and the physical and topological properties of the Schwarz and Gyroid surfaces. We use the similar methodology as in \cite{Asb221}, constructing several supervised learning models to learn the relations between several numerical characteristics obtained from the surfaces. This methodology allows a mathematician to quickly <<test>> the intuition about the possible relationship between two quantities and the way they may be related.

During our experiments, we were able to construct 4 models (2 models for each surface) that provide a clear relationship between 
\begin{itemize}
    \item shape factor and porosity,
    \item shape factor and 1-persistence entropy.
\end{itemize}

To accomplish this, we utilized the {\it Python} {\it Scikit-Learn} \cite{scikit-learn} library, a built-in Python module for machine learning. 
To split the data into a test and training set and to increase the reliability of the required models, we employed the \textup{cross\_val\_score} and KFold functions from the \textup{model\_selection} module for further hypothesis testing. 
The dataset was partitioned into training and testing subsets using KFold cross-validation \cite{scikit-learn}. The employment of KFold cross-validation addresses the necessity of testing models by extracting the test data from the dataset. Random selection of test subsets results in significant variations in RMSE data, depending on the dataset split. To mitigate this randomness, KFold cross-validation was applied with \texttt{n\_splits=4}, as the total number of points in the dataset is divisible by 4.

As a quality metric, we will use the Root Mean Squared Error (RMSE), which is defined as:

\[RMSE = \sqrt{\frac{1}{n}\sum_{i=1}^{n}(y_i - \hat{y}_i)^2}\]

Where \(n\) represents the number of observations, \(y_i\) is the actual value for the i-th observation, and \(hat{y}_i\) is the predicted value for the i-th observation.

\subsection{Shape factor predicts persistence}
Let's examine the relationship between the shape factor and persistence. We have noticed that the shape factor crucially changes the surface and its shape, therefore we hypothesised that there's an undiscovered relationship between the shape factor of a surface and its topology. To numerically test the hypothesis, we vectorized the persistence diagrams obtained in Section \ref{subsec:ph_of_tpms}, as discussed in Section \ref{subsec:entropy}. We have tested several vectorization methods, but in the end we decided to use persistence entropy as the vectorization method. 

We utilized the polynomial regression model as a supervised learning model. Since we only employ shape factor as an input feature for the model, our model simply boils down to a univariate polynomial of some degree. Since the degree of the model is a hyperparameter, we conducted the experiments with the models of different degrees. 

Besides the degree of the model, since persistence entropy is computed only in a specific dimension, we need to choose the dimension. We were unable to construct any model that would detect the possible relation between the shape factor and persistence entropy in either $0$-th or $2$-nd dimension, so we've focused only on the peristence entropy in dimension $1$. Table \ref{tab:entropy} presents the 1-persistence entropy values for both Schwarz and Gyroid models with respect to the shape factor.

\begin{longtable}{|c|c|c|c|}
    \caption{Persistence entropy calculation results}
    \label{tab:entropy} \\
    \hline
        № & Shape factor & $\operatorname{PE}_1$ of Schwarz & $\operatorname{PE}_1$ of Gyroid\\ \hline 
        \endfirsthead
        
        \multicolumn{4}{c}%
        {\tablename\ \thetable\ -- \textit{Continued from previous page}} \\
        \hline
        1 & 2 & 3 & 4\\ \hline
        \endhead
        
        \hline
        \multicolumn{4}{r}{\textit{Continued on next page}} \\
        \endfoot
        
        \hline
        \endlastfoot
        
1 & -0.985086 & 7.828019 & 12.817731 \\ \hline
2 & -0.834337 & 8.086322 & 11.421037 \\ \hline
3 & -0.602925 & 8.421331 & 10.059915 \\ \hline
4 & -0.357099 & 8.694553 & 8.853321 \\ \hline
5 & -0.350885 & 8.698730 & 8.856987 \\ \hline
6 & -0.262321 & 8.778409 & 8.563515 \\ \hline
7 & -0.086151 & 8.893371 & 8.089089 \\ \hline
8 & 0.184648 & 8.821489 & 7.165093 \\ \hline
9 & 0.429486 & 8.609666 & 6.492101 \\ \hline
10 & 0.446178 & 8.594761 & 6.458588 \\ \hline
11 & 0.465395 & 8.577064 & 6.398883 \\ \hline
12 & 0.583936 & 8.452268 & 6.092819 \\ \hline
13 & 0.642550 & 8.380563 & 5.921346 \\ \hline
14 & 0.713229 & 8.281657 & 5.739262 \\ \hline
15 & 0.721315 & 8.267385 & 5.713075 \\ \hline
16 & 0.778426 & 8.178185 & 5.532658 \\ \hline
17 & 0.867778 & 8.030560 & 5.232656 \\ \hline
18 & 0.873779 & 8.022236 & 5.235386 \\ \hline
19 & 0.953246 & 7.872525 & 5.054533 \\ \hline
20 & 0.959561 & 7.865700 & 5.054643 \\ \hline
\end{longtable}

\begin{figure}[H]
  \centering
    \includegraphics[width=0.85\textwidth]{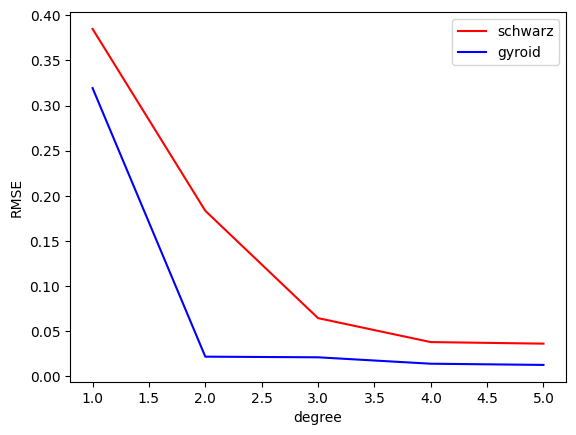}
    \caption{Root Mean Square Error (RMSE) for varying polynomial degrees.}
    \label{fig:entropy_rmse}
\end{figure}

Table \ref{tab:polynomial_model_entropy} displays the polynomial regression models' performance with different degrees on the test subset, and Figure \ref{fig:entropy_rmse} illustrates the variation in RMSE on test data with changing polynomial degrees. This represents the mean RMSE values obtained from different splits in KFold cross-validation.

According to Table \ref{tab:polynomial_model_entropy}, the minimum mean test RMSE for the Schwarz surface is achieved with a polynomial degree of $4$, and for the Gyroid model, with a degree of $2$.

Another approach is to train polynomial models on the entire dataset, calculate the RMSE, and employ the elbow method to determine the degree value beyond which RMSE reduction is primarily due to overfitting. Figure \ref{fig:entropy_rmse} depicts this graphically. Visually, the optimal degree value for the Gyroid surface is $2$, and for the Schwarz model, it is $3$ or $4$.

This method is inherently subjective and unreliable. Nonetheless, both the elbow method and the preceding calculations provide strong evidence that, at least for the Gyroid surface, a degree $2$ polynomial can describe the relationship between the shape factor and the 1-persistence entropy.

\begin{longtable}{|c|c|c|c|}
    \caption{Polynomial regression model performance for shape factor and persistence entropy}
    \label{tab:polynomial_model_entropy} \\
    \hline
    № & Degree & Mean Test RMSE for Schwarz & Mean Test RMSE for Gyroid \\ \hline
    \endfirsthead
    \multicolumn{4}{c}%
        {\tablename\ \thetable\ -- \textit{Continued from previous page}} \\
    \hline
    № & Degree & Mean Test RMSE for Schwarz & Mean Test RMSE for Gyroid \\ \hline
    \endhead
    \hline
    \multicolumn{4}{r}{\textit{Continued on next page}} \\
    \endfoot
    \endlastfoot
    1 & 1 & 0.410427 & 0.418837 \\ \hline
    2 & 2 & 0.255528 & 0.025773 \\ \hline
    3 & 3 & 0.117023 & 0.055177 \\ \hline
    4 & 4 & 0.062917 & 0.029832 \\ \hline
    5 & 5 & 0.068506 & 0.133599 \\ \hline
    6 & 6 & 0.220459 & 0.272051 \\ \hline
    7 & 7 & 0.395311 & 0.226815 \\ \hline
    8 & 8 & 0.910539 & 1.387404 \\ \hline
\end{longtable}

Beyond the polynomial model, we have also experimented with an exponential model. Such model was supposed to detect a following relationship hypothesis:

$$ \text{There exists a polynomial } P \text{ and a constant } c \text{ such that } \operatorname{PE}_1(M) = ae^{P(d(M))}. $$

Table \ref{tab:exponential_model_entropy} presents the mean RMSE values across KFold splits on test data with varying polynomial degrees. According to the table, the minimum mean test RMSE for the Schwarz surface is achieved with a polynomial degree of $3$, and for the Gyroid surface, with a degree of $2$. The degree value for the Gyroid surface aligns with the polynomial model, with a slight difference for the Schwarz surface.

\begin{longtable}{|c|c|c|c|}
    \caption{Exponential regression model performance for shape factor and persistence entropy}
    \label{tab:exponential_model_entropy} \\
    \hline
    № & Degree & Mean Test RMSE for Schwarz & Mean Test RMSE for Gyroid \\ \hline
    \endfirsthead
    \multicolumn{4}{c}%
        {\tablename\ \thetable\ -- \textit{Continued from previous page}} \\
    \hline
    № & Degree & Mean Test RMSE for Schwarz & Mean Test RMSE for Gyroid \\ \hline
    \endhead
    \hline
    \multicolumn{4}{r}{\textit{Continued on next page}} \\
    \endfoot
    \endlastfoot
    1 & 1 & 0.242320 & 0.424048 \\ \hline
    2 & 2 & 0.305091 & 0.021064 \\ \hline
    3 & 3 & 0.074146 & 0.040753 \\ \hline
    4 & 4 & 0.112339 & 0.028396 \\ \hline
    5 & 5 & 0.379184 & 0.124916 \\ \hline
    6 & 6 & 5.164287 & 0.269900 \\ \hline
    7 & 7 & 0.498869 & 0.192177 \\ \hline
    8 & 8 & 1.340072 & 0.690750 \\ \hline
\end{longtable}

Notably, there is little difference between polynomial and exponential models in the minimal mean RMSE for both surfaces. Since the exponential model essentially builds upon the polynomial one, we discard the exponential model hypothesis. 

Thus, these experimental results provide strong evidence of a potential polynomial relationship between the shape factor and the 1-persistence entropy for the Gyroid and Schwarz surfaces. In particular, we obtained the following relationships in the conducted experiments:

\begin{align*}
    &\text{Schwarz:} &&\operatorname{PE}_1 = 7.74 - 2.92d + 0.52d^2 - 1.06d^3 + 0.66d^4,\\
    &\text{Gyroid: } &&\operatorname{PE}_1 = 8.83 - 1.06d^2,\\
\end{align*}
where $d$ is the shape factor and $\operatorname{PE_1}$ is the 1-persistence entropy. These considerations lead us to the following conjecture:

\begin{conjecture}
    For any triply-periodic minimal surface $M$, there exists a polynomial $P \in \mathbb{R}[t]$ such that 
    $$ \operatorname{PE}_1(M) = P(d(M)), $$
    where $\operatorname{PE}_1(M)$ is the 1-persistence entropy and $d(M)$ is the shape factor.
\end{conjecture}

\subsection{Shape factor predicts porosity}

Similar to the point described above, let us analyze the relationship between the aspect ratio and porosity in the Schwarz and Gyroid TPMS. Our hypothesis is that there is an undetected relationship between these measures. We put forward the following hypothesis:

$ \text{There exists some polynomial } P \text{ such that } Pr(M) = P(d(M)), $

\textit{where} \(Pr(M)\) is TPMS porosity function from object M, \(P(d(M))\) - this is some polynomial of the function of the shape factor and the object M.

As in the example described above, using cross-validation, we divided the training and test set data into 4 groups and trained and tested models for different groups. The table \ref{tab:polynomial_model_porosity} presents the results of supervised learning the Schwarz and Gyroid TPMS.

\begin{longtable}{|c|c|c|c|}
    \caption{Polynomial regression model for shape factor and porosity}
    \label{tab:polynomial_model_porosity} \\
    \hline
        № & Degree & Mean test RMSE for Schwarz & Mean test MSE for Gyroid\\ \hline 
        \endfirsthead
        
        \multicolumn{4}{c}%
        {\tablename\ \thetable\ -- \textit{Continued from previous page}} \\
        \hline
        1 & 2 & 3 & 4\\ \hline
        \endhead
        
        \hline
        \multicolumn{4}{r}{\textit{Continued on next page}} \\
        \endfoot
        
        \hline
        \endlastfoot
        
1 & 1 & 0.035765 & 0.049244 \\ \hline
2 & 2 & 0.015448 & 0.021227 \\ \hline
3 & 3 & 0.015445 & 0.154869 \\ \hline
4 & 4 & 0.210800 & 0.249551 \\ \hline
5 & 5 & 0.079919 & 0.990721 \\ \hline
6 & 6 & 2.419949 & 2.262439 \\ \hline
7 & 7 & 18.010812 & 6.147779 \\ \hline
8 & 8 & 121.490580 & 74.947167 \\ \hline
\end{longtable}

\begin{figure}[H]
  \centering
    \includegraphics[width=0.85\textwidth]{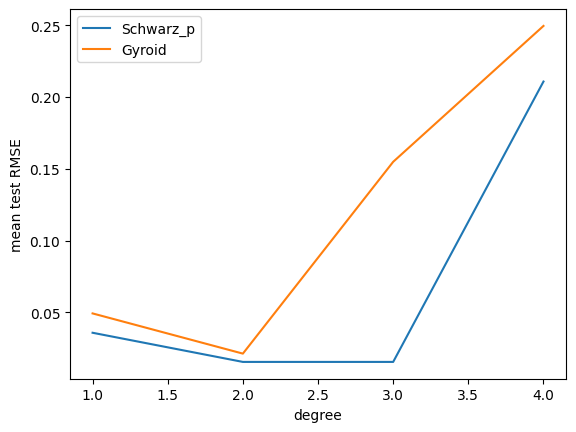}
    \caption{The mean test RMSE values between shape factor and porosity per degree of the polynomial for two type of surfaces}
    \label{fig:polynomial_model_porosity}
\end{figure}

Table \ref{tab:polynomial_model_porosity} presents the quality of the polynomial regression models with varying degrees on the test subset of shape factor and porosity. The figure \ref{fig:polynomial_model_porosity} depicts the change in RMSE on the test data when modifying the degree. This represents the mean RMSE value obtained from different splits in K-Fold cross-validation.

Therefore, the best Gyroid models have a second degree polynomial and have the following mathematical form:

\[Pr_{gyroid}(M) = 0.09674384d(M) - 0.06645933d(M)^2\]

For the Schwarz surface, models with the 2nd and 3rd degree of polynomial showed better quality, but the model with the third degree of polynomial turned out to be slightly better. Therefore, their mathematical representation is:

\[Pr_{schwarz}(M) = -0.06816784d(M) + 0.04984266d(M)^2 + 0.01013945d(M)^3\]
and 
\[Pr_{schwarz}(M) = -0.06150581d(M) + 0.05036231d(M)^2\]

\begin{figure}[H]
  \centering
    \includegraphics[width=0.85\textwidth]{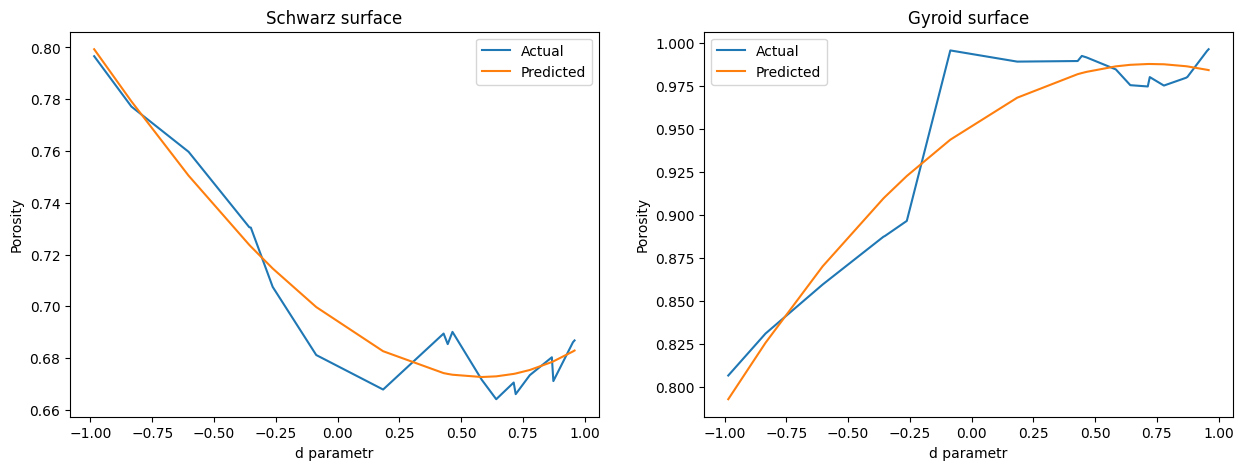}
    \caption{Comparison of the best polynomial regression models of Schwarz and Giroyd surfaces with their real data on the dependence of the shape factor (d parametr) and porosity}
    \label{fig:compare_with_real_data}
\end{figure}

Let's plot the results of the most accurate polynomial regression models and compare them with real data. This comparison can be seen in the graphs \ref{fig:compare_with_real_data}.

 A hypothesis was proposed about the existence of a polynomial dependence between the shape and porosity of the Schwarz and Gyroid surfaces. The results of cross-validation showed that high-quality polynomial regression models can be constructed for both types of surfaces, confirming the proposed hypothesis. The obtained models demonstrate good agreement with the real data, which indicates the presence of the identified relationship between the shape and porosity of the triply periodic minimal surfaces.

%% file: conclusion.tex
\section{Conclusion}\label{106}

Triply Periodic Minimal Surfaces (TPMS) have attracted significant attention due to their unique properties. These surfaces are characterized by a high surface area-to-volume ratio and mathematically controlled geometry, finding widespread applications in medical implants \cite{Asb28}, various engineering fields \cite{Asb29}, as well as energy-absorbing materials \cite{Asb27} and other industrial applications \cite{Asb211, Asb212, Asb210}. The design of TPMS-based materials and structures involves careful topological optimization to achieve target physical properties corresponding to specific application requirements.

We conducted a analysis of both the physical and topological properties of TPMS models, depending on their mathematical representation. The obtained results can be considered as a promising direction in the tasks of computer modeling and design of objects requiring the properties inherent in TPMS.

During the study, we have discovered two relationships between the shape factor and both 1-persistence entropy and porosity. We have constructed the polynomial regression models for the Schwarz and Gyroid surfaces, that model such relationships with a sufficiently small values of the RMSE quality metric. The results are presented in Table \ref{tab:conclusion_results}, where the porosity of the triply periodic minimal surfaces serves as the predictor.

\begin{longtable}{|p{0.1cm}|p{1.5cm}|p{2cm}|p{5.5cm}|c|}
    \caption{Results of dependence models}
    \label{tab:conclusion_results} \\
    \hline
    № &  Surface &  Predictor &  Mathematical model &  RMSE\\ \hline 
        \endfirsthead
        
        \multicolumn{5}{c}%
        {\tablename\ \thetable\ -- \textit{Continued from previous page}} \\
        \hline
        1 & 2 & 3 & 4 & 5 \\ \hline
        \endhead
        
        \hline
        \multicolumn{5}{r}{\textit{Continued on next page}} \\
        \endfoot
        
        \endlastfoot
        
        1 & Schwarz & porosity & $Pr_{schwarz}(M) = 0.694 - 0.068d(M) + 0.05d(M)^2 + 0.01d(M)^3$ & 0.015 \\ \hline
        2 & Schwarz & 1-persistence entropy & $\operatorname{PE_1}(M) = 7.74 - 2.92d(M) + 0.52d(M)^2 - 1.06d(M)^3 + 0.66d(M)^4$ & $0.063$ \\ \hline
        3 & Gyroid & porosity & $Pr_{gyroid}(M) = 0.953 + 0.097d(M) - 0.067d(M)^2$ & 0.021 \\ \hline
        4 & Gyroid & 1-persistence entropy & $\operatorname{PE_1}(M) = 8.83 - 1.06d(M)^2$ & $0.026$ \\ \hline
\end{longtable}

The developed models are planned to be used in further TPMS-related research. Their application will allow the design of specific forms and structures with the desired properties inherent in minimal surfaces, which, in turn, will contribute to the solution of various applied problems. Additionally, further research should be aimed at improving the obtained models, as well as exploring new ways of applying deep learning and machine learning methods to achieve even better results.

On another hand, the conducted research provides an example of the use of machine learning techniques in mathematics in the spirit of \cite{Asb221}. We have hypothesised that there's an undiscovered relationship, and supervised models were able to detect such relationships, at least in the case of aforementioned surfaces. As far as we know, this work is the first one that propose such interconnections.

In future research we plan to prove the constructed conjectures, thus mathematically establishing the discovered relationships. Another fruitful direction in the future work is to expand this methodology in order to uncover another patterns in the topology of TPMS.



